\documentclass[a4paper,12pt,aps,superscriptaddress]{article}
\usepackage{amsmath}

\usepackage{amsfonts}
\usepackage{amssymb}

\usepackage{float}
\usepackage{epsfig}

\newcommand{\be}{\begin{equation}}
\newcommand{\ee}{\end{equation}}
\newcommand{\bea}{\begin{eqnarray}}
\newcommand{\eea}{\end{eqnarray}}

\begin{document}

\title{{\bf Statistical Curse of the  Second Half Rank, \\ Eulerian numbers\\and Stirling numbers}
$$ $$
St\'ephane Ouvry\\
 Universit\'e  Paris-Sud, Laboratoire de Physique Th\'eorique et Mod\`eles
Statistiques\\  
  UMR 8626, F-91405 Orsay}

\maketitle

\begin{abstract}
I describe the occurence of Eulerian numbers and  Stirling  numbers of the second kind in the combinatorics of the "Statistical Curse of the  Second Half Rank" problem \cite{0}. 

\medskip\noindent {}

\end{abstract}

\section{Introduction}

The "Statistical Curse of the  Second Half Rank" problem \cite{0} stems  from real life considerations leading to  rather complex combinatorics. One is primarly concerned  with
 rank expectations in  sailing boats regattas, bearing in mind that the issue discussed here is quite general and can apply  to rank expectations of students taking exams, or other types of similar endeavors as well. 
 \begin{figure}[htbp]
 \begin{center}
\epsfxsize=8cm
\centerline{\epsfbox{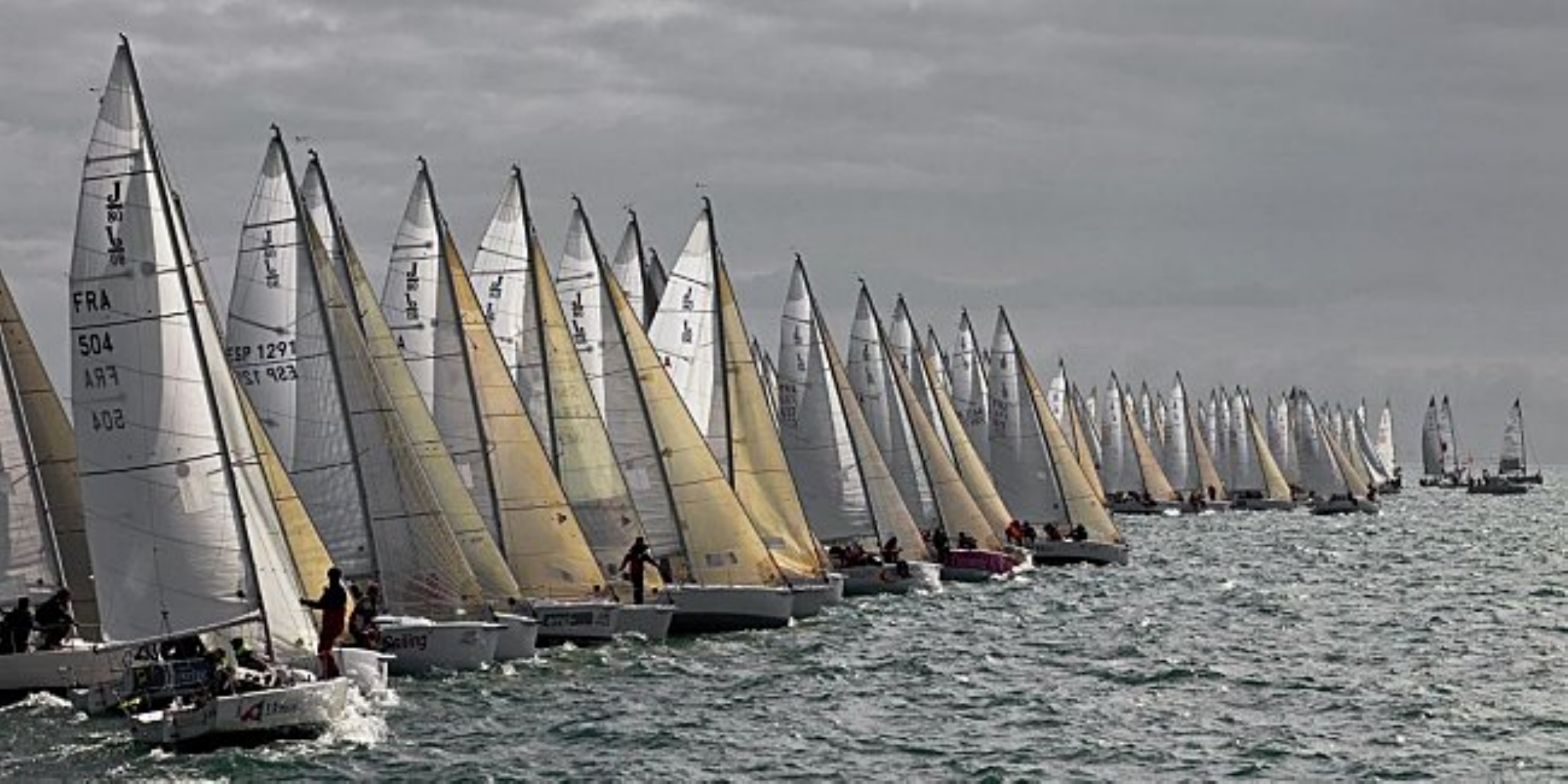}}
\end{center}
\caption{Just before the start of a race at the {\it Spi Ouest-France}: a big number of identical boats is going to cross a virtual starting line in a few seconds.}
 \end{figure}

 Consider as a typical example
  the {\it Spi Ouest-France} regatta  which takes place each year during  4 days at Easter  in La Trinit\'e-sur-Mer,   Brittany (France).    
   It involves a "large" number $n_b$ of  identical boats, say    $n_b= 90$, 
running a "large" number $n_r$ of races, say $n_r= 10 $ (that is to say  $2,3$ races  per day, weather permitting, see Fig.~1).
In each race each boat gets a rank $ 1\le {\rm rank}\le 90$ with the condition that there are
no ex aequo.
Once the last race is over, to determine   the final rank of a boat and thus  the winner of the regatta one proceeds as follows:

 1)   one adds  each boat's rank in each race    $\rightarrow$ its score $n_t$: 
 here $n_b=90, n_r=10$ so that $10\le  n_t\le 900$
 
\noindent  $n_t=10$ is the  lowest possible score $\to$ the boat was  always ranked $1^{\rm rst}$  

\noindent  $n_t=900$ is the highest possible score $\to$ the boat was  always ranked  $90^{\;\rm th}$

\noindent $n_t=10\times (1+90)/2\simeq 450$ the middle score $\to$   the boat was     on average ranked $45^{\;\rm th}$ 

 2)   one orders the boats according to  their score $\rightarrow$ their final rank
 
\noindent the boat with the lowest score  $\to$ $1^{\rm rst}$ (the   winner)  
 
\noindent the boat with  next to the lowest score  $\to$ $2^{\rm nd}$
 
 \noindent etc...
  
  \begin{figure}[htbp]
\begin{center}
\epsfxsize=16cm
\centerline{\epsfbox{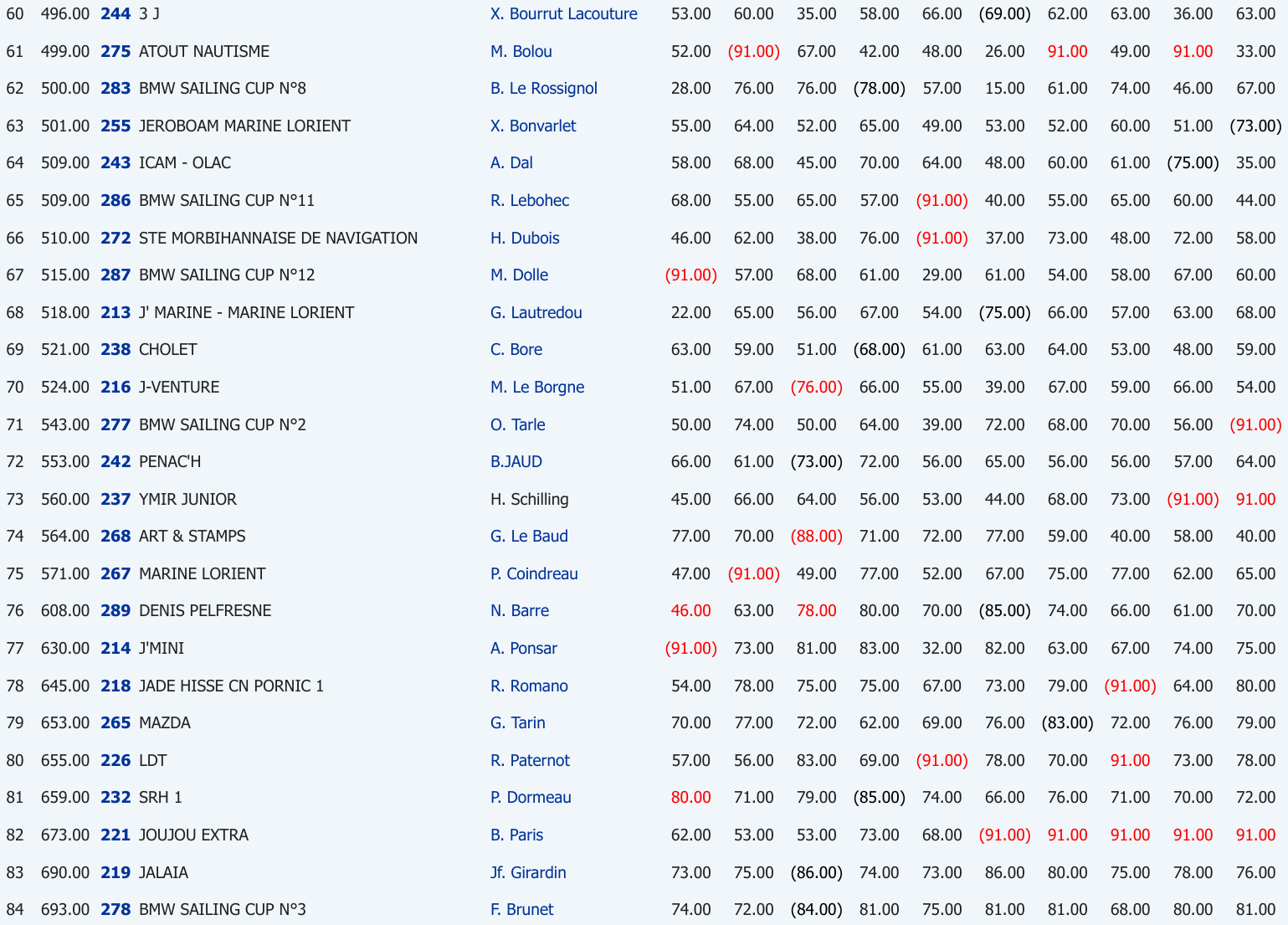}}
\end{center}
 \caption{Starting from  the left the first  column gives the final rank of the boat, the second column its score,  the third column its "improved" score,  the fourth column its name, the fifth column the name of its skipper, and the  next $10$ columns  its  ranks in each of the $10$ races.}
\end{figure}

  What is the "Statistical Curse of the  Second Half Rank"?    
  In  the {\it Spi Ouest-France}  2009  results sheet  (see Fig.~2 for  boats with final rank between $60^{\rm th}$  and  $84^{\rm th}$), consider for example      the boat with final rank $70^{\rm th}$:
 its  ranks in the $10$ races are $51$, $67$, $76$, $66$, $55$, $39$, $67$, $59$, $66$, $54$ so that its score is $ n_t=600$. The crew of this boat might   naively expect that since
  its  mean rank is ${600/ 10}=60$ and so
  it has been on average ranked $60^{\rm th}$,  its final rank should be around $60^{\rm th}$.
No way, the boat ends\footnote{This pattern  is even more pronounced if one  notes in Fig.~2  that  the rank  in the third race $(76)$   is  in parenthesis: this is because   each boat's worst  rank  is removed from the final counting. It is as if they were
 only  $9$ races with, for the boat considered here, ranks $51, 67,  66, 55, 39, 67, 59, 66, 54$, "improved"  score   $n_t=524$  and   mean rank ${524/ 9}\simeq 58$. So, on average, the boat is ranked $58^{\rm th}$ even though it ends up being $70^{\rm th}$.}  up being $70^{\rm th}$. 
 
This "curse" phenomenon is quite general and one would like to understand its origin   and be able to evaluate it.
A qualitative explanation  is  simple \cite{0}: in a given race
given the  rank of the boat considered above, 
 assume   that the ranks of the other boats  are random variables with a uniform distribution. 
The  
 random rank assumption is   good  if, bearing in mind that all boats are identical,  the crews can also be considered as more or less equally worthy, which for sure is partially the case.
Since there are no ex aequo it means that the
 ranks of the other boats,  in the first race, are  a random permutation of $(1, 2, 3,\ldots, 50, 52,\ldots, 90)$, 
 in the second race, a random permutation of $(1, 2, 3,\ldots, 66, 68,\ldots, 90)$, etc.

 Each  race is obviously independent from  the others,  so that 
 the  scores  are the sums of $10$ independent random variables.
But $10$ is already a large number in probability calculus so that 
the Central Limit Theorem  applies. It follows that the  scores  are random variables  with a gaussian probability density centered around the middle score $\simeq 450$.
 A gaussian distribution  implies  a lot of boats with scores packed around the middle score. 
   Since the score $600$ of the boat considered here  is larger than the middle score $450$, this packing implies that its final rank is  pushed upward from its mean rank: {\bf this is the statistical "curse"}. On the contrary if the boat's score  had been lower than the middle score, its final rank would have been pushed downward from its mean rank : {\bf its crew would have enjoyed  a statistical "blessing"}.

Let us 
  rewrite things  more precisely by asking, 
 given the score $n_t$ of the boat considered among the $n_b$ boats,
 what is the probability distribution $P_{n_t}(m)$ for its final rank    to be $m\in[1,n_b]$?

A complication  arises as soon as $n_r\ge 3$:  
$P_{n_t}(m)$ does not only  depend  on the score  $n_t$   but also on the ranks of the boat in each race.  
For example for $n_r= 3$,  take $n_b=3$  and    the score $n_t=6$, 
it is  easy to check by complete enumeration that  
$P_{6=2+2+2}(m)\ne P_{6=1+2+3}(m) $. The   distributions are of course similar but slightly differ. 
To avoid this complication let us  from now on
consider  $n_b$ boats with  in each race  random ranks given by  a  random permutation  of $(1,2,3,\dots,n_b)$ 
$\oplus$ an additional "virtual" boat  only specified by its score $n_t$ and ask  
the question again: given the score $n_t$ of this virtual boat what is   the probability distribution $P_{n_t}(m)$ for its final rank to be  $m\in[1,n_b+1]$ ? 
This problem is  almost the same\footnote{In the $2$-race case it is in fact the same problem.}  yet a little bit simpler since, by construction, it does not have the complication discussed above.

 In a given race $k$ call $n_{i,k}$ the rank of the boat $ i$   with $1\le i\le n_b$  and $1\le k\le n_r$.
There are no ex aequo in a given  race:  the $n_{i,k}$'s are a random permutation of $(1,2,3,\dots,n_b)$ so that   they  are correlated random variables with
 
$${\rm sum\; rule}\quad \quad{\sum_{i=1}^{n_b}n_{i,k}}=1+2+3+\ldots+n_b={n_b(1+n_b)\over 2}$$ 
 
$${\rm mean} \quad \quad\langle n_{i,k}\rangle={1+n_b\over 2}$$ 
 
$${\rm fluctuations}\quad \quad \langle n_{i,k}n_{j,k}\rangle-\langle n_{i,k}\rangle\langle n_{j,k}\rangle={1+n_b\over 12}(n_b\delta_{i,j}-1)$$

\noindent Now, the score $ n_i$ of  boat $i$ is defined as $ n_i\equiv \sum_{k=1}^{n_r}n_{i,k}$, the  middle score being  $n_r(1+n_b)/ 2 $.
In the large $n_r$ limit the Central Limit Theorem applies -here for correlated random variables-
to yield the scores joint density probability  distribution

$$f(n_1,\dots,n_{n_b})=$$   

$$
\sqrt{2\pi \lambda n_b } \left(\sqrt{1\over 2\pi\lambda}\right)^{n_b}\delta\left(\sum_{i=1}^{n_b} (n_i-n_r{1+n_b\over 2})\right)\exp\left[
-{1\over 2\lambda}\sum_{i=1}^{n_b}(n_i-n_r{1+n_b\over 2})^2\right]$$
 with $\lambda=n_r {n_b(1+n_b)/ 12}$. 

 Now consider the virtual boat with  score $n_{t}$: 
 $P_{n_t}(m)$  is  
the probability for $m-1$ boats among the $n_b$'s to have a  score $n_i<n_{t}$  and for the other $n_b-m+1$ boats to have a score $n_i\ge n_{t}$
 $$
P_{n_{t}}(m)={n_b\choose m-1}\int_{-\infty}^{n_{t}}dn_1 \dots dn_{m-1}
\int_{n_{t}}^{\infty}dn_{m} \dots dn_{n_b}f(n_1,\dots ,n_{n_b})$$

\noindent Let us  take the large   number of boats limit:  
a saddle point approximation finally \cite{0} gives   
 $\langle m\rangle$ as the cumulative probability distribution of a normal
variable
$$
 {\langle m\rangle}= \frac{n_b}{\sqrt{2\pi\lambda}} \int_{-\infty}^{\bar{n}_{t}} \exp\left[ -\frac{n^2}{2\lambda}
\right] dn
$$
where
$$\bar{n}_{t}= n_{t}-n_r\frac{(1+n_b)}{2} $$
and  $ n_r\le n_t\le n_r n_b\rightarrow  - n_r{n_b\over 2}\le\bar{n}_{t}\le n_r{n_b\over 2}$. In Fig.~3 a plot of ${\langle m\rangle/n_b}$ is displayed in the case 
$n_r=30$, $n_b=200$. The curse and blessing effects are clearly visible on the sharp increase around the middle score -a  naive expectation would  claim a linear increase.
\begin{figure}[htbp]
\begin{center}
\epsfxsize=8cm
\centerline{\epsfbox{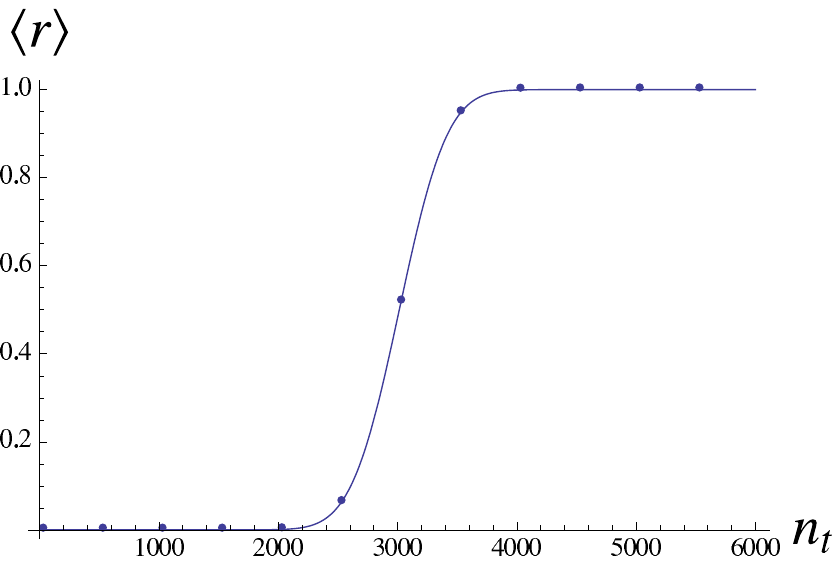}}
\end{center}
\caption{$\langle r\rangle = {\displaystyle \langle m\rangle\over\displaystyle  n_b}$, $n_r=30$,
$n_b=200$, 
middle score $3000$,
dots = numerics}
\end{figure}
The variance  can be obtained along similar lines with a manifest damping  due to the correlations effects.

So far we have dealt with the "curse" which is a large number of races and boats effect. 
Let us now turn to the combinatorics  of a small number of races $n_r=2,3,\ldots$ for a  given number of boats $n_b=1,2,\ldots$.
The simplest situation is the $2$-race case $n_r=2$ which happens to be solvable -it can be viewed as a solvable "$2$-body" problem- with an exact solution  for  $P_{n_t}(m)$. 
To obtain in this simple situation the probability distribution  one proceeds as follows  (see Fig.~4):  

1) one represents the possible  ranks configurations of any boat among the $n_b$'s  in the two races by points on a $n_b\times n_b$ square lattice 
(in the $3$-race  case one would have a cubic lattice, etc):
 since there are no ex aequo there is exactly  $1$ point per line and per column, so there are $n_b!$ such configurations
(in the $n_r$ races case one would have   $(n_b!)^{n_r-1}$ such configurations)

2) one enumerates all the configurations with $m-1$ points below  the  diagonal $n_t$:
this is the number of configurations with final rank $m$ for the virtual boat\footnote{Contrary to the no ex aequo rule in  a given race, boats of course have equal final ranks if they have the same score. }.

\begin{figure}[htbp]
\begin{center}
\epsfxsize=8cm
\centerline{\epsfbox{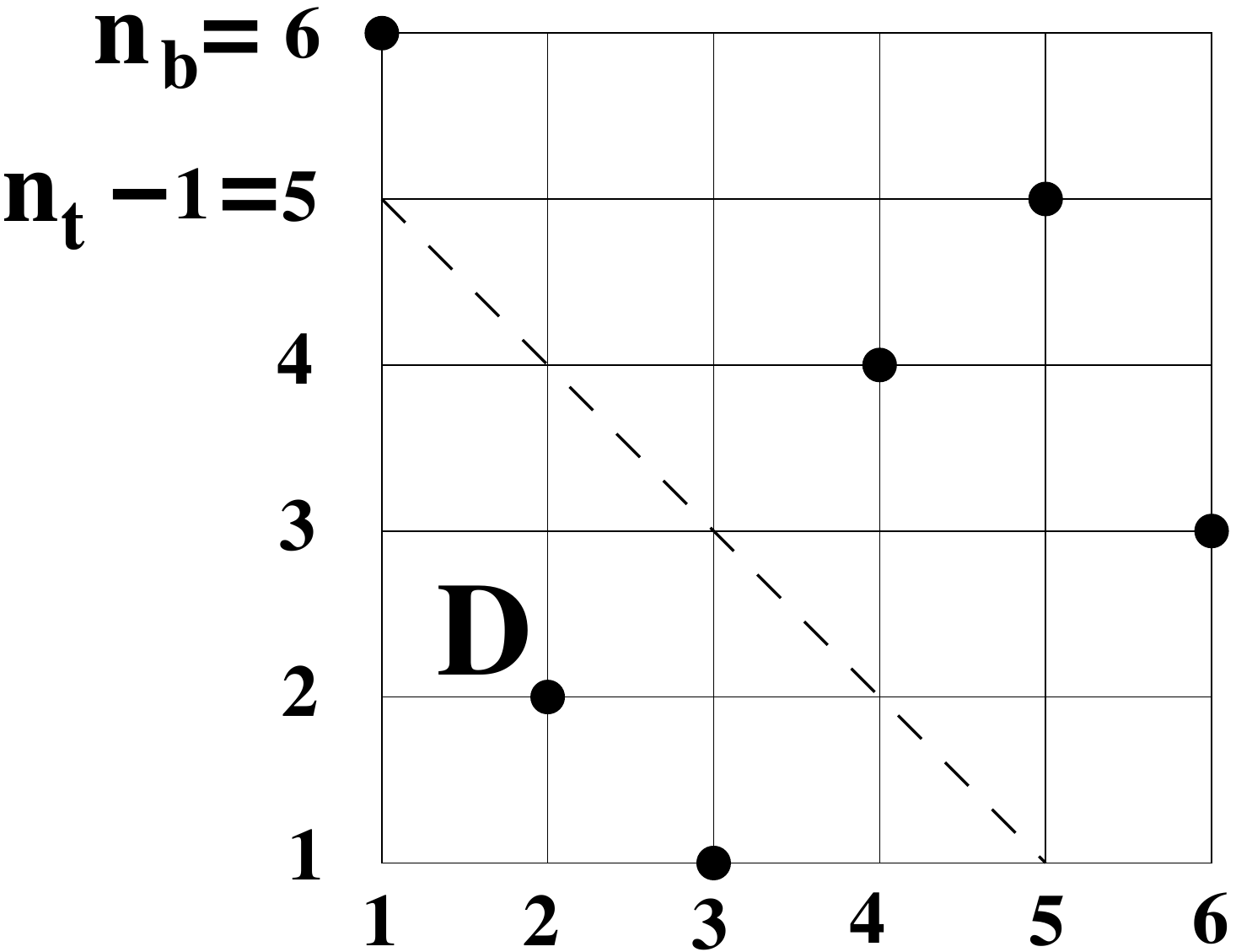}}

\end{center}
\caption{ The $2$-race case: a $m=3$ configuration for $n_b=6$ and $n_t=6$. The dashed line is the diagonal $n_t=6$.}
 
\end{figure}

The problem has been narrowed down to a combinatorial enumeration which is  doable  \cite{0}: 
one finds for $2\le n_t\le 1+n_b$ 

$$
P_{n_t}(m)=(1+n_b)\sum_{k=0}^{m-1}(-1)^k(1+n_b-n_t+m-k)^{n_t-1}\frac{(n_b-n_t+m-k)!}{
 k! (1+n_b-k)! (m-k-1)!    }         
$$

\noindent and for $ n_t=n_b+1+i\in[n_b+2, 2n_b+1]$ with $i=1,2,\ldots, n_b$, by symmetry, $P_{n_b+1+i}(m)=P_{n_b+2-i}(n_b +2-m)$.

\vspace{0.4cm}
So far no particular numbers, be they Eulerian or Stirling, have occured.  
Let us concentrate on  the  middle score $n_t=2{(1+n_b)/ 2}=1+n_b$ to get

$$ P_{n_t=1+n_b}(m)=(1+n_b)\sum_{k=0}^{m-1}  (-1)^k  \frac{(m-k)^{n_b}}{ k! (1+n_b-k)! }  
$$

Let us tabulate  $ P_{n_t=1+n_b}(m)$, with $m\in[1,n_b+1]$, for $n_b=1,2,...7$:
$$ \{1,0\}$$
$$ {1\over 2!}\{1,1,0\}$$
$$ {1\over 3!}\{1,4,1,0\}$$
$$ {1\over 4!}\{1,11,11,1,0\}$$
$$ {1\over 5!}\{1,26,66,26,1,0\}$$
$$ {1\over 6!}\{1,57,302,302,57,1,0\}$$
$$ {1\over 7!}\{1,120,1191,2416,120,1,0\}$$
The numbers between  brackets happen to be known as  the   ${\rm Eulerian}(n_b,k)$ numbers with $k=m-1\in[0,n_b-1]$ (here one has dropped the trivial $0$'s obtained for $m=n_b+1$ i.e. $k=n_b$). An Eulerian number (see Fig.~5) is 
the number of permutations of the numbers 1 to n in which exactly m elements are greater than the previous element (permutations with m "ascents") as illustrated  in Fig.~6.
\begin{figure}[htbp]
\begin{center}
\epsfxsize=8cm
\centerline{\epsfbox{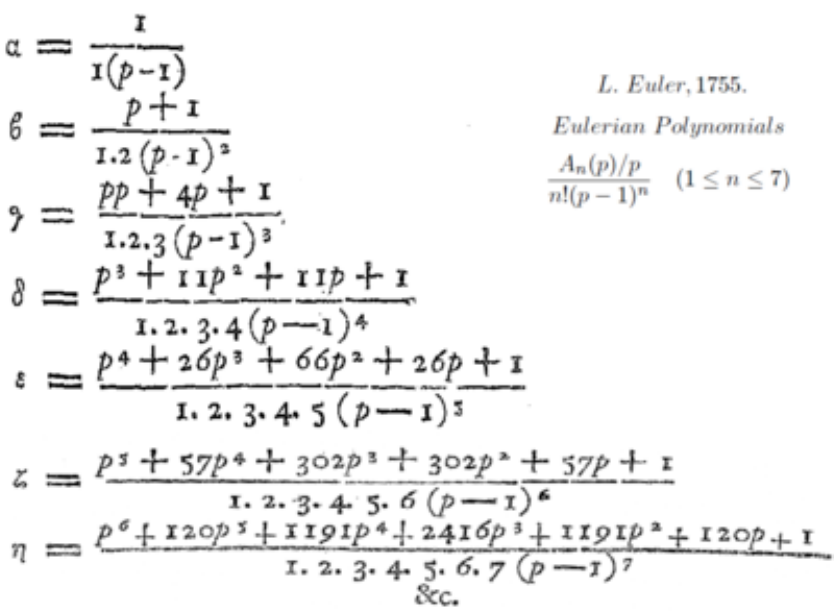}}
\end{center}
\caption{The Eulerian numbers.}
\end{figure}
\begin{figure}[htbp]
\begin{center}
\epsfxsize=8cm
\centerline{\epsfbox{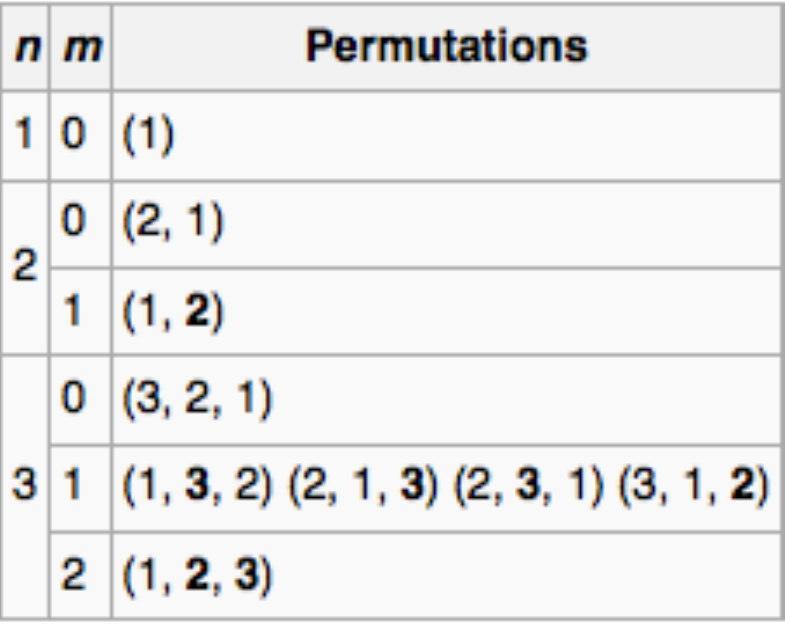}}
\end{center}
\caption{ Ascents for $n=1,2,3$. For $n=4$,  the permutation $(1, 4, 2, 3)$   has $ m=2$ ascents.}
\end{figure}
Their generating function is
$$ g(x,y)=\frac{e^x-e^{x y}}{e^{x y}-e^x y}$$
with a series expansion  $x\simeq 0$
$$ g(x,y)\simeq x+\frac{1}{2!} x^2 (y+1)+\frac{1}{3!} x^3 \left(y^2+4
   y+1\right)+\frac{1}{4!} x^4 \left(y^3+11 y^2+11
   y+1\right)$$
   $$+\frac{1}{5!} x^5 \left(y^4+26 y^3+66 y^2+26
   y+1\right)+\frac{1}{6!} x^6 \left(y^5+57 y^4+302 y^3+302
   y^2+57 y+1\right)+O\left(x^7\right)$$

It is not difficult to realize that all the information in $P_{n_t}(m)$ is  contained in $P_{n_t=1+n_b}(m)$ that is to say in the generating function $g(x,y)$. Rephrased more precisely,  $P_{n_t=n_b+1}(m)$  for $n_b=1,2,\ldots$  is generated by  $y g(x,y)$  with the $x$ exponent being the number of boats $n_b$  and, for a given $n_b$,  the $y$ exponent being the rank $m\in[1,n_b+1]$. Similarly
 $P_{n_t=n_b}(m)$ for $n_b=2,3,\ldots$ is  generated by

$$(y-1) \int_0^x g(\text{z},y) \, d\text{z}+g(x,y)-x =\frac{e^x-e^{x y}}{e^{x y}-e^x y}+\frac{(y-1) \left(\log
   (1-y)-\log \left(e^{x y}-e^x y\right)\right)}{y}-x$$
   $$=x^2 y+\frac{1}{3} x^3 y (y+2)+\frac{1}{12} x^4 y \left(y^2+7
   y+4\right)+\frac{1}{60} x^5 y \left(y^3+18 y^2+33
   y+8\right)$$
   $$+\frac{1}{360} x^6 y \left(y^4+41 y^3+171
   y^2+131 y+16\right)+O\left(x^7\right) $$
This procedure can be repeated for $n_t=n_b-1, n_b-2, \ldots$ with expressions involving double, triple, $\ldots$ integrals of $g(x,y)$.

Why Eulerian numbers should play a role  here
can be  explained from scratch by a simple combinatorial argument\footnote{I thank S. Wagner for drawing my attention to this explanation.}: let us use the notations that   the boat which came up $i^{\rm th}$ in the first race had rank  $a(i)$ in the second race: $(a(1),a(2),\ldots, a(i),\ldots, a(n_b))$ is a random permutation of $(1, 2,\ldots, i,\ldots, n_b)$. In the case  of interest where the score of the virtual boat is the middle score $n_b+1$, the boat that came $i^{\rm th}$ in the first race beats the virtual  boat if $a(i)\le n_b-i$, that is to say if it is better than $i^{\rm th}$, counting from the bottom, in the second race. Therefore, the counting problem of $P_{n_t=n_b+1}(m)$ is equivalent to counting so-called excedances: an excedance in a permutation $(a(1),a(2),\ldots, a(i),\ldots, a(n_b))$ is an element such that $a(i) > i$. The number of permutations with precisely $m$ excedances is known to be an Eulerian number (thus excedances are what is called an Eulerian statistic, see for example \cite{wag} p.~23). As an illustration look at the case  $n_b=3$: the six permutations of $(1,2,3)$ are $(1,2,3)$, $(1,3,2)$, $(2,1,3)$, $(2,3,1)$, $(3,1,2)$, $(3,2,1)$. The numbers of excedances (here defined as  $a(i)\le n_b-i$) are respectively $1,1,2,1,1,0$ which indeed yields the ${\rm Eulerian}(3,k)$  numbers  
$$1,4,1$$ 
  for $k=0,1,2$ ($1$ permutation with 0 excedance, $4$ permutations with $1$ excedance, $1$ permutation with $2$ excedances) 
 appearing  in $P_{n_t=1+n_b}(m)$, $m\in[1,n_b+1]$, $n_b=3$.

 There is still an other way   to  look at  the problem in terms of  
 Stirling numbers of the second kind, here  defined as
$$ n_{n_t}(i)=\frac{1}{(-i+{n_t}-1)!}\sum _{j=1}^{{n_t}-1} j^{{n_t}-2}
   (-1)^{-i-j+{n_t}}
   \binom{-i+{n_t}-1}{j-1}$$
 Stirling numbers count in how many ways  the numbers $(1,2,\ldots,n_t-1)$ can be partitioned in $i$ groups: for example for $n_t=5$ 

\hspace{2.4cm}$\to 1$  way  to split the numbers $(1,2,3,4)$ into $4$ groups  

$(1),(2),(3),(4)$

\hspace{2.4cm}$\to 6$ ways to split the numbers $(1,2,3,4)$ into $3$ groups

$(1),(2),(3,4);(1),(3),(2,4);(1),(4),(2,3);(2),(3),(1,4);(2),(4),(1,3);(3),(4),(1,2)$

\hspace{2.4cm}$\to 7$ ways to split the numbers $(1,2,3,4)$ into $2$ groups
 
$(1),(2,3,4);(2),(1,3,4);(3),(1,2,4);(4),(1,2,3);(1,2),(3,4);(1,3),(2,4);(1,4),(2,3)$

\hspace{2.4cm}$\to 1$  way  to split the numbers $(1,2,3,4)$ into $1$ group

$(1,2,3,4)$

\noindent so that one obtains  

$${1,6,7,1}$$
 The probability distribution $P_{n_t}(m)$ can indeed be 
 rewritten  \cite{duduche}  as 
$$P_{n_t}(m)={1\over n_b!}\sum_{i=m}^{n_t-1}(-1)^{i+m}{n_{n_t}(i)}(1+n_b-i)!\bigg({i-1\atop m-1}\bigg)$$
Why Stirling  numbers should play a role  here
arises \cite{duduche} 
 from graph counting considerations 
on the  $n_b\times n_b$ lattice  when one now includes all the points below the diagonal.
 As an example  let us still consider the case $n_t=5$:  below the diagonal $n_t=5$ they are $6$  points $a,b,c,d,e,f$ labelled by their lattice coordinates  $a=(1,1),\; b=(1,2),\; c=(1,3),\; d=(2,1),\; e=(2,2),\; f=(3,1)$.  One  draws a graph according to the  no ex aequo exclusion rule: starting say from the point $a$  one links it to another point if it  obeys  the  no ex aequo exclusion rule with respect to $a$, that is to say if its coordinates are not $(1,1)$.  In our example there is only one such point $e$. Then one can link the point $e$ to the  points $c$ and $f$, which can also be linked together. Finally one can link $c$  to $d$ and $f$ to $b$ with also a  link between $d$ and $b$. The numbers ${n_{n_t=5}(i+1)}$ count in the graph just obtained  the number of subgraphs with either $i=1$ point (this is the number of points $6$), $i=2$ points linked (there are 7 such cases), $i=3$ points fully linked (there is 1 such case), $i=4$  points fully linked (there is no such case), $\ldots$, a counting which  finally gives the Stirling-like numbers 
 $$1,6,7,1$$ 
 where the $1$ on the left is by convention the number of subgraphs with $i=0$ point.


It still remains to be shown why this subgraph counting is indeed equivalent to the Stirling counting. To do  this one has simply to notice that the former  is encapsulated in a  recurrence relation obtained by  partitioning the $n_{n_t+1}(i+1)$  lattice counting  as:

either there is $0$  point on the diagonal $n_t$ $\to n_{n_t}(i+1)\left({n_t-1\atop 0}\right)$  

either there is $1$  point on the diagonal $n_t$ $\to n_{n_t-1}(i)\left({n_t-1\atop 1}\right)$ 

either there are $2$  points on the diagonal $n_t$ $\to n_{n_t-2}(i-1)\left({n_t-1\atop 2}\right)$ 

etc

\noindent It follows that the numbers $ n_{n_t+1}(i+1)$ have to obey  the recurrence relation
 $$   n_{n_t+1}(i+1)=\sum_{k'=0}^{i} n_{n_t-k'}(i+1-k')\left({n_t-1\atop k'}\right)   $$
\noindent valid for $i\in[0,n_t-2]$, bearing in mind that when $i=n_t-1$ trivially  $ n_{n_t+1}(n_t)=1$.
 But this recurrence relation can be  mapped on a more standard recurrence relation for the Stirling numbers of the second kind: 
if one sets $k=n_t-i\in[2,n_t]$  and defines  now  ${\rm Stirling}(n_t,k)\equiv n_{n_t+1}(i+1)$ one obtains for the ${\rm Stirling}(n_t,k)$'s
 $$   {\rm Stirling}(n_t,k)=\sum_{k'=0}^{n_t-k} {\rm Stirling}(n_t-k'-1,k-1)\left({n_t-1\atop k'}\right)   $$
 which via $k''=n_t-k'\in[k,n_t]$ rewrites as 
 $$   {\rm Stirling}(n_t,k)=\sum_{k''=k}^{n_t} {\rm Stirling}(k''-1,k-1)\left({n_t-1\atop k''-1}\right)   $$
that is to say finally
$$   {\rm Stirling}(n_t+1,k+1)=\sum_{k''=k}^{n_t} {\rm Stirling}(k'',k)\left({n_t\atop k''}\right)   $$

\noindent This therefore establishes that the ${\rm Stirling}(n_t,k)$'s, i.e.  the ${n_{n_t}(i)}$'s, are indeed  Stirling numbers of the second kind.
This is of course not a surprise that both Eulerian and Stirling numbers do play a role since there is a  correspondance\footnote{For  Eulerian, Stirling, and other well-known numbers in combinatorics see for example \cite{wag}.} between them
$$ {\rm Eulerian}(n_b,k)=\sum _{j=1}^{k+1}(-1)^{k-j+1}\left({n-j\atop n-k-1}\right)j!\;{\rm Stirling}(n_b,j)$$
with $k\in[0,n_b-1]$.

Why  does one  rewrite  $P_{n_t}(m)$  in terms of Stirling numbers in the 2-race case? Because this rewriting can be  generalized \cite{duduche}   to the $n_r$-race case. The generalisation is formal since one does not know what  the     $n_r$-dependant "generalized Stirling" numbers which  control the probability distribution\footnote{The lattice becomes larger with the number of races: as said above when $n_r=3$ one has a cubic lattice.  The corresponding graph structure becomes more an more involved, and so the counting of the associated fully connected subgraphs.} are. Again this is like moving from a solvable "$2$-body" problem to a so far non solvable "$n_r$-body" problem.

Still it remains quite fascinating that well-known numbers in combinatorics, such as Eulerian and Stirling numbers, should be at the heart of the understanding of rank expectations in  regattas, at least in the 2-race case.

{\bf Acknowledgements:} I would like to take the occasion of  Leonid Pastur $75^{\rm th}$ birthday celebration to tell  him all my friendship and admiration. I acknowledge  useful conversations with S. Wagner. My thanks also to Dhriti Bhatt's youthful energy and her willingness to look  at this problem again in the summer 2012.

\end{document}